\documentclass[letterpaper, 10 pt, conference]{ieeeconf}  

\IEEEoverridecommandlockouts                              
\usepackage{cite, graphicx, amsmath, amssymb, mathrsfs, float, subfig, color, url}
\newtheorem{theorem}{Theorem}
\renewcommand{\baselinestretch}{1.0}

\title{\LARGE \bf
Robust Control Framework for Time-Varying Power-Sharing among Distributed Energy Resources
}

\author{Mayank Baranwal$^{1,a}$ and Srinivasa M. Salapaka$^{1,b}$
\thanks{$^{1}$Department of Mechanical Science and Engineering, University of Illinois at Urbana-Champaign, 61801 IL, USA}
\thanks{$^{a}${\tt\small baranwa2@illinois.edu}, $^{b}${\tt\small salapaka@illinois.edu}}
}

\begin{document}

\maketitle
\thispagestyle{empty}
\pagestyle{empty}

\begin{abstract}
One of the most important challenges facing an electric grid is to incorporate renewables and distributed energy resources (DERs) to the grid. Because of the associated uncertainties in power generations and peak power demands, opportunities for improving the functioning and reliability of the grid lie in the design of an efficient, yet pragmatic distributed control framework with guaranteed robustness margins. This paper addresses the problem of output voltage regulation for multiple DC-DC converters connected to a grid, and prescribes a robust scheme for sharing power among different sources. More precisely, we develop a control architecture where, unlike most standard control frameworks, the desired power ratios appear as reference signals to individual converter systems, and not as internal parameters of the system of parallel converters. This makes the proposed approach suited for scenarios when the desired power ratios vary rapidly with time. Additionally, the proposed control framework is suitable to both centralized and decentralized implementations, i.e., the same control architecture can be employed for voltage regulation irrespective of the availability of common load-current (or power) measurement, without the need to modify controller parameters. The control design is obtained using robust optimal-control framework. Case studies presented show the enhanced performance of prescribed optimal controllers for voltage regulation and power sharing.
\end{abstract}

\section{INTRODUCTION}\label{sec:Intro}
\begin{figure}[!t]
	\centering
	\includegraphics[width=0.95\columnwidth]{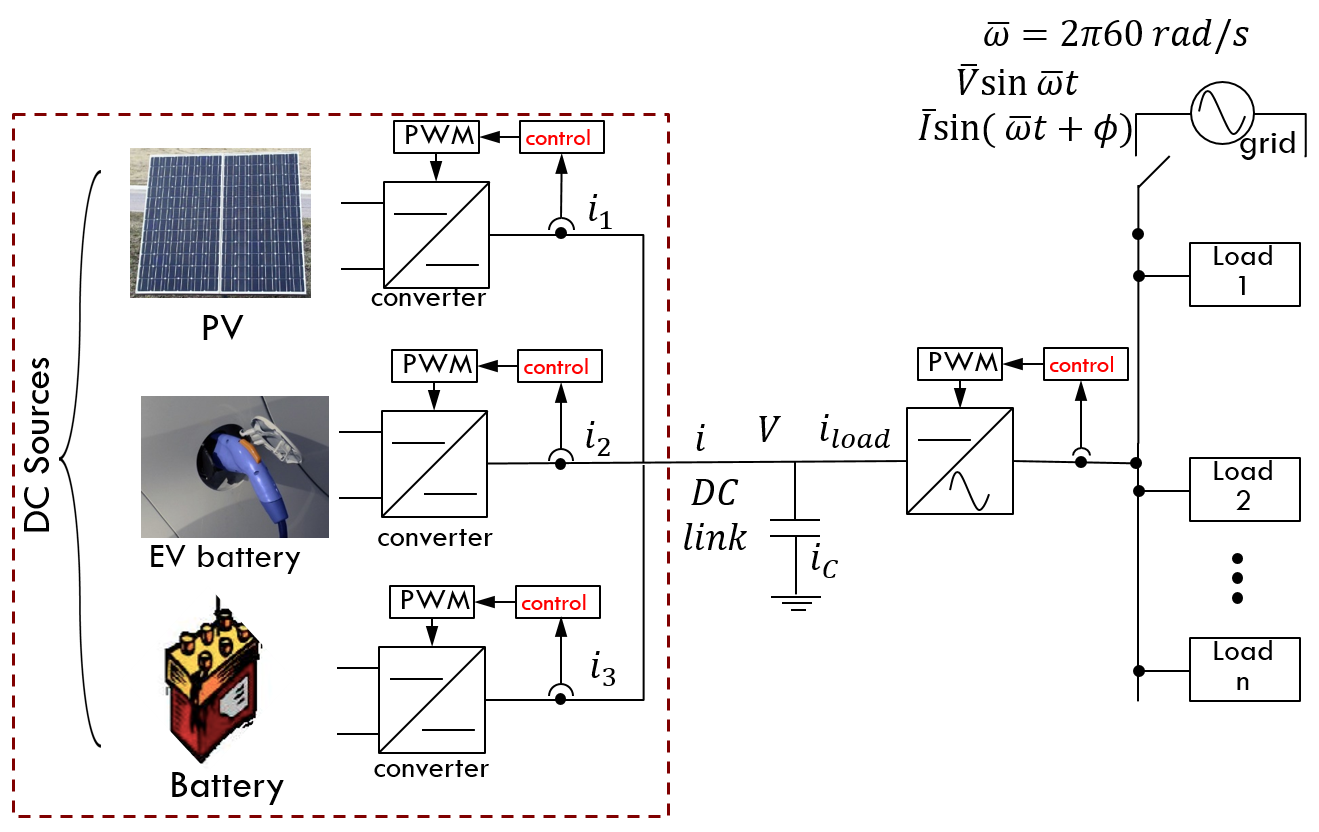}
	\vspace{-0.5em}
	\caption{{\small A schematic of a microgrid. An array of DC sources provide power for AC loads. Power sources  provide  power at DC-link, their common output bus, at a voltage that is regulated to a set-point. The control system at the respective DC-DC converter that interfaces with a source is responsible for regulating the voltage at he DC-link. An inverter that connects to the DC-link converts the total current from the sources at the regulated voltage to  alternating current (AC) at its output to satisfy the power demands of the AC loads. This paper describes an approach for control design of the multiple converters systems associated with power transfer from sources to the DC-link (shown by the dotted line).}} 
	\label{fig:MCS}
	\vspace{-1.5em}
\end{figure}

High environmental impact of fossil-based energy sources and demand for future energy sustainability have resulted in increased interest in the use of renewable energy resources such as solar and wind. Integration of renewable resources on the community level is achieved through smart microgrids. Microgrids are localized grid systems that are capable of operating in parallel with, or independently from, the existing traditional grid \cite{lasseter2002microgrids, hatziargyriou2007microgrids}. Microgrid technology enables integration of renewable energy sources such as solar and wind energy, distributed energy resources (DERs), energy storage, and demand response. Fig. \ref{fig:MCS} shows a schematic of a microgrid with multiple DC sources providing power for AC loads. Although in most traditional grids that rely on conventional sources of {\em dispatchable} electric power, the power output of renewables can not be manipulated. Limited predictability with such resources result in intermittent power generation; moreover time-varying loads, practicability and economics factors pose additional challenges in efficient operation of microgrids. Thus it is required to develop efficient distributed control technologies for reliable operation of smart microgrids.

In such smart grids, multiple DC power sources connected in parallel, each interfaced with DC-DC converter, provide power at their common output, the DC-link, at a regulated voltage; this power can directly feed DC loads or be used by an inverter to interface with AC loads. In this paper, we address the problem of distributed control of DC-DC converters for output voltage regulation, and time-varying power sharing (dictated by the economic layer) among multiple power sources. The main challenges arise from the uncertainties in the size and the schedules of loads, the complexity of a coupled multi-converter network, the uncertainties in the model parameters at each converter, and the adverse effects of interfacing DC power sources with AC loads, such as the $120$ Hz ripple that has to be provided by the DC sources.

Problems pertaining to robust and optimal control of converters have received recent attention. Conventional PID-based controllers often to fail address the problem of robustness and modeling uncertainties. In \cite{olalla2010lmi}, a linear-matrix-inequality (LMI) based based robust control design for boost converters has resulted in significant improvements over PID based control designs. In \cite{weiss2004h, hornik2011current, salapaka2014viability}, robust $\mathcal{H}_\infty$-control framework is employed in the context of inverter systems. While the issue of current sharing is extensively studied \cite{panov1997analysis, wu1993load}, most methods assume a single power source. A systematic control design that addresses all the challenges and objectives for the multi-converter control is still lacking. The control architecture proposed in this paper addresses the following primary objectives - 1) voltage regulation at the DC-link with guaranteed robustness margins, 2) prescribed time-varying power sharing in a network of parallel converters, 3) controlling the trade-off between $120$Hz ripple on the total current provided by the power sources and the ripple on the DC-link voltage. While these objectives are partially addressed in our prior work \cite{baranwal2016robust} on the robust control of DC-DC converters, a main drawback of the design proposed in \cite{baranwal2016robust} is that the control framework does not allow for time-varying power sharing requirements. In this work, we propose a novel control framework wherein the power requirements on each converter are imposed through external references, and thus the framework allows for time-varying power sharing by incorporating high-bandwidth robust controllers.

The control architecture proposed in this work exploits structural features of the paralleled multi-converter system, which results in a modular and yet coordinated control design. For instance, noting that the objective of voltage regulation is common to all converters; accordingly at each converter, it employs a nested (outer-voltage inner-current) control structure \cite{erickson2007fundamentals}, where all converters share the same design for the outer-loop voltage controllers while the inner-loop current controllers are so chosen that the entire closed-loop multi-converter system can be reduced to an equivalent single-converter system in terms of the transfer function from the desired regulation setpoint $V_{ref}$ to the voltage at the DC-link $V_{dc}$. The controllers are designed for fully {\em centralized} implementation with the instantaneous load current $i_{load}$ measurement accessible to all converters; however, in practice $i_{load}$ is often {\em estimated} through power calculations on the AC side and communicated to individual converters only at a rate slower than the sampling rate of the controllers. In the wake of this limitation, we propose a novel method for voltage regulation and power sharing that is inspired by conventional voltage droop method. An interesting aspect of the proposed implementation is that the same outer controllers $K_{v}$ and $K_{r}$ along with the shaped inner plant $\tilde{G}_{c,n}$ can be employed even for the scenario where $i_{load}$ measurement is unavailable. Thus the proposed framework is applicable to both {\em centralized} and {\em decentralized} implementations. An important revelation provided by the application of $\mathcal{H}_\infty$ robust optimal control is the underlying optimal structure in the outer-loop controllers $K_{v}$ and $K_{r}$. While we are yet to explore the reasoning behind the hidden optimal structure, it helps in further reduction of the overall complexity of the distributed control design from analysis and implementation points of view.

The rest of the paper is organized as follows. Sec. \ref{sec:model} describes the averaged modeling of DC-DC converters. We then describe the control design methodology for a single converter system in Sec. \ref{sec:Control_Single}, followed by its extension to a network of parallel converters in Sec. \ref{sec:Control_Many}. The underlying theory is then corroborated by extensive simulations in Sec. \ref{sec:simulations}, followed by some important conclusions and immediate directions to future works.

\section{MODELING OF CONVERTERS}\label{sec:model}
\begin{figure}
	\subfloat[]{\includegraphics[width= 0.9\columnwidth]{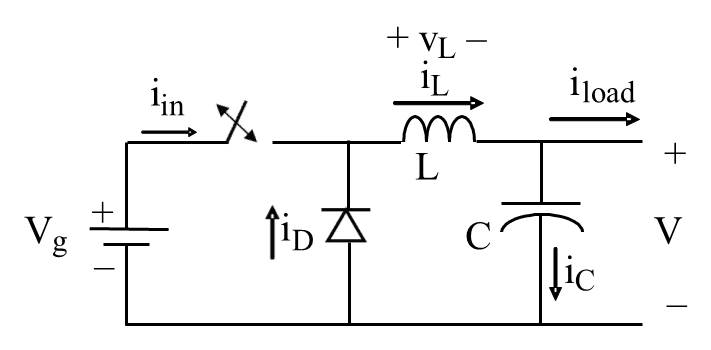}}\\
	\vspace{-.5em}
	\subfloat[]{\includegraphics[width= 0.9\columnwidth]{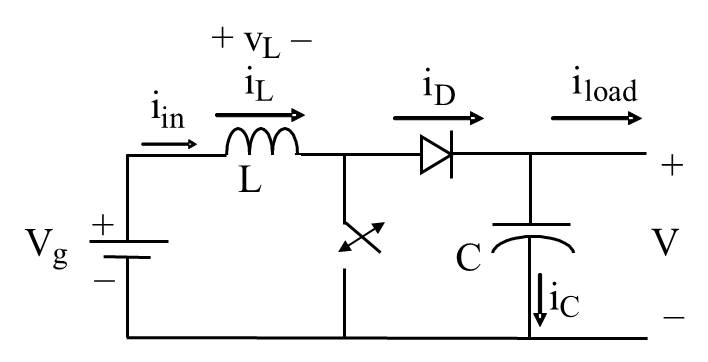}}
	\caption{{\small Schematics for (a) Buck converter, (b) Boost Converter. The converters are assumed to operate in continuous-conduction-mode (CCM).}}
	\vspace{-1.5em}
	\label{fig:conv_schematic}
\end{figure}
In this section, we describe the differential equations that govern the dynamics of DC-DC converters. These converters belong to a class of switched-mode power electronics, where a semiconductor based high-frequency switching mechanism (and associated electronic circuit) connected to a DC power source enables changing voltage and current characteristics at its output. The models presented below depict dynamics for signals that are averaged over a switch cycle.

Fig. \ref{fig:conv_schematic}a shows a schematic of a Buck converter. Buck converter regulates a voltage at its output which is lower than the input voltage. The averaged dynamic model of a Buck converter is given by
\begin{small}
	\begin{eqnarray}\label{eq:buck_eq}
		L\frac{di_L(t)}{dt} &=& \underbrace{-V(t)+d(t)V_g}_{\tilde{u}(t):= -V(t)+u(t)} = \tilde{u}(t)\nonumber\\
		C\frac{dV(t)}{dt} &=& i_L(t) - i_{load}(t),
	\end{eqnarray}
\end{small}
where $d(t)$ represents the duty-cycle (or the proportion of {\em ON} duration) at time $t$. Note that the prescribed averaged model does not explicitly require any information on the output load.

Similarly we can describe the averaged dynamics of a Boost converter (shown in Fig. \ref{fig:conv_schematic}b), given by
\begin{small}
	\begin{eqnarray}\label{eq:boost_eq}
		L\frac{di_L(t)}{dt} = \underbrace{V_g-d'(t)V(t)}_{\tilde{u}(t):= V_g-u(t)} = \tilde{u}(t)&& \nonumber\\
		C\frac{dV(t)}{dt} = \underbrace{\left(D'+\hat{d'}\right)}_{\approx D'}i_L(t) - i_{load}(t) \approx D'i_L(t) - i_{load}(t),&&
	\end{eqnarray}
\end{small}
where $d'(t):=1-d(t)$ and $D'=\left(V_g/V_{ref}\right)$. Here $V_{ref}$ represents the desired output voltage. Note that $\hat{d}(t) = d'(t)-D'$ is typically very small, and therefore allows for a linear approximation around the nominal duty-cycle, $D = 1-D'$. We use both $V$ and $V_{dc}$ interchangeably to denote the output voltage at the DC-link.

\section{PROBLEM FORMULATION}\label{sec:problem}
This paper addresses the following primary objectives {\em simultaneously} - (1) Output voltage regulation in presence of time-varying loads/generation and parametric uncertainties, (2) time-varying current (power) sharing among multiple sources, (3) $120 Hz$ ripple current sharing between inductor current $i_L$ and capacitor current $i_C$. The last two objective is dealt in our prior work \cite{baranwal2016robust} and is addressed by an appropriate design of inner-controller described in Sec. \ref{subsec:Single_inner} and \ref{sec:Control_Many}. In this paper, we primarily focus on achieving the first two objectives, while inheriting the properties of the inner-controller for ripple current sharing.

\section{CONTROL FRAMEWORK FOR SINGLE CONVERTER}\label{sec:Control_Single}
\begin{figure}
	\includegraphics[width= 0.9\columnwidth]{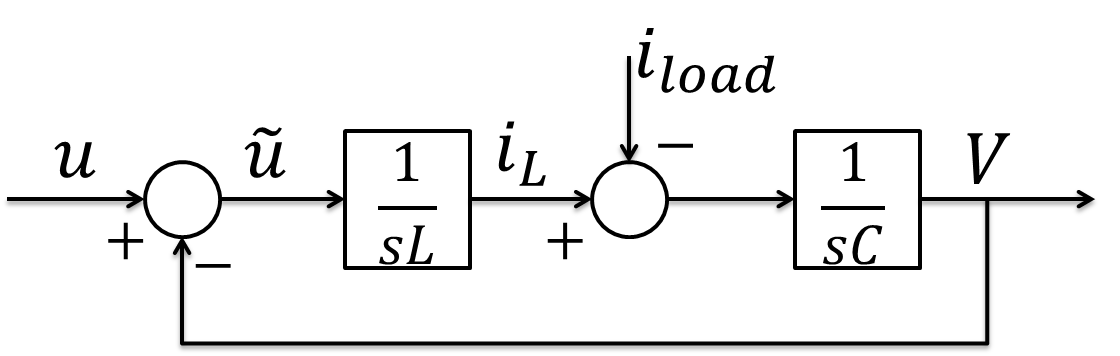}
	\caption{{\small Block diagram representation of a buck-type converter. The control signal $\tilde{u}$ is converted to an equivalent PWM signal to command the gate of the transistor acting as a switch.}}
	\vspace{-1.5em}
	\label{fig:buck_model}
\end{figure}

In this section, we describe the {\em inner-outer} controller architecture for a single Buck converter system. The corresponding block diagram representation of the dynamical equations in \ref{eq:buck_eq} is shown in Fig. \ref{fig:buck_model}. While the design is easily extendable to include other converter types such as Boost and Buck-Boost, the discussion has been confined to Buck converters only for the sake of brevity. Unlike most other {\em outer}-controllers, the {\em outer}-controller proposed in this paper takes into account both DC-link voltage $V$ and load current $i_{load}$ measurements. The requirements on current sharing are imposed through this additional $i_{load}$ measurement (as explained in Sec. \ref{sec:Control_Many}).

\subsection{Design of the inner-loop controller}\label{subsec:Single_inner}
The design for the inner-loop controller $K_c$ is inherited from our previous work \cite{baranwal2016robust}. The inner-loop controller $K_c$ is designed as a $2^{nd}$-order controller to ensure a relatively low-order optimal controller $K_v = \left[K_{v},K_{r}\right]^T$ (see Fig. \ref{fig:single_block}). The main objective for designing the inner-loop controller $K_c$ is to decide the trade-off between the $120$Hz ripple on the capacitor current $i_C$ (equivalently on the output voltage $V_{dc}$) and the inductor current $i_L$ of the converter. Accordingly, $K_c$ is designed such that the inner-shaped plant $\tilde{G}_c$ is given by
\begin{small}
	\begin{equation}\label{eq:Gc_tilde}
		\tilde{G}_c(s) = \left(\frac{\tilde{\omega}}{s+\tilde{\omega}}\right)\left(\frac{s^2+2\zeta_1\omega_0s+\omega_0^2}{s^2+2\zeta_2\omega_0s+\omega_0^2}\right).
	\end{equation}
\end{small} 
where $\omega_0 = 2\pi 120$rad/s and $\tilde{\omega}, \zeta_1, \zeta_2$ are design parameters. The parameter $\tilde{\omega} > \omega_0$ and it is used to implement a low-pass filter to attenuate undesirable frequency content in $i_L$ beyond $\tilde{\omega}$. Thus, the bandwidth of the inner-shpaed plant is decided by the choice of $\tilde{\omega}$. The parameters $\zeta_1$ and $\zeta_2$ impart a notch-like behavior to $\tilde{G}_c$ at $\omega_0 = 120$Hz, and the size of the notch is determined by the ratio $\zeta_1/\zeta_2$. Note that $\tilde{G}_c$ represents the inner closed-loop plant from the output of the outer-loop controllers $u$ to the inductor current $i_L$, and since $i_C = i_L - i_{load}$, the ratio $\zeta_1/\zeta_2$ can be appropriately designed to achieve a specified trade-off between $120$Hz ripple on $i_C$ and $i_L$. The stabilizing $2^{nd}$-order controller $K_c$ that yields the aforementioned inner closed-loop plant $\tilde{G}_c$ is explicitly given by
\begin{small}
	\begin{equation}\label{eq:Kc}
		K_c(s) = L\tilde{\omega}\frac{(s^2+2\zeta_1\omega_0s+\omega_0^2)}{(s^2+2\zeta_2\omega_0s+\omega_0^2+2(\zeta_2-\zeta_1)\omega_0\tilde{\omega})}.
	\end{equation}
\end{small}
The readers are encouraged to refer to Sec. III in \cite{baranwal2016robust} for further details on the inner-loop control design.

\subsection{Design of the outer-loop controller}\label{subsec:Single_outer}
\begin{figure}
	\includegraphics[width= 0.95\columnwidth]{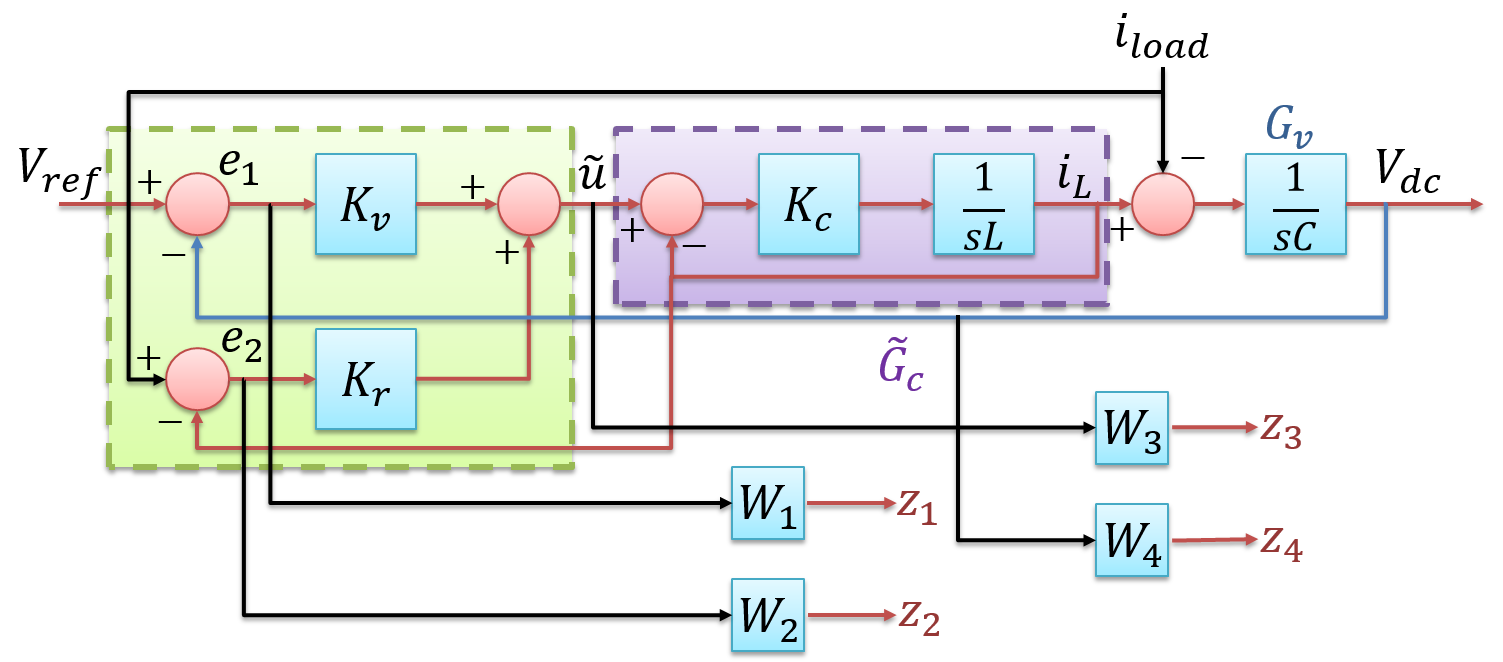}
	\caption{{\small Block diagram representation of the inner-outer control design. Exogenous signal $V_{ref}$ represents the desired output voltage. The quantities $V_{dc}$, $i_{load}$ and $i_L$ represent the available measurements.}}
	\vspace{-1.5em}
	\label{fig:single_block}
\end{figure}
For a given choice of inner-controller $K_c$, we present our analysis and design of controller in terms of transfer function block diagrams shown in Fig. \ref{fig:single_block}. In this figure, $\tilde{G}_c$ represents the inner shaped plant. The outer controllers are denoted by $K_{v}$ and $K_{r}$, and are designed to regulate the output DC voltage $V_{dc}$ to the desired reference voltage $V_{ref}$ and the inductor current $i_L$ to the load current $i_{load}$, respectively. Note that from (\ref{eq:buck_eq}), $i_L$ is equal to $i_{load}$ at steady-state. The augmentation of controller $K_{r}$ forms the basis for time-varying power sharing and is explained in the next section.

The performance of a Buck converter is characterized by its voltage and power reference tracking bandwidths, better voltage signal-to-noise-ratio (SNR), and robustness to modeling uncertainties. The main objective for the design of the controllers $K_{v}$ and $K_{r}$ is to make the tracking errors small and simultaneously attenuate measurement noise to achieve high resolution. This is achieved by posing a {\em model-based} multi-objective optimization framework, where the required objectives are described in terms of norms of the corresponding transfer functions, as described below. From Fig. \ref{fig:single_block}, the transfer function from exogenous inputs $w = \left[V_{ref}, i_{load}, \tilde{u}\right]^T$ to regulated output $z = \left[z_1, z_2, z_3, z_4, e_1, e_2\right]$ is given by
\begin{small}
	\begin{equation}\label{eq:hinf}
	\left[ \begin{array}{c}
		z_1\\
		z_2\\
		z_3\\
		z_4\\
		e_1\\
		e_2
	\end{array} \right] = 
	\left[ \begin{array}{ccc}
		W_1 & W_1G_v & -W_1G_v\tilde{G}_c \\
		0 & W_2 & -W_2\tilde{G}_c \\
		0 & 0 & W_3 \\
		0 & -W_4G_v & W_4G_v\tilde{G}_c \\
		1 & G_v & -G_v\tilde{G}_c \\
		0 & 1 & -\tilde{G}_c
	\end{array}	\right] 
	\left[ \begin{array}{c}
		V_{ref}\\
		i_{load}\\
		\tilde{u}
	\end{array} \right].
	\end{equation}
\end{small}
The optimization problem is to find stabilizing controllers $K_{outer} = \left[K_{v},K_{r}\right]^T$ such that the $\mathcal{H}_\infty$-norm of the above transfer function from $w$ to $z$ is minimized. Here the weights $W_1, W_2, W_3$ and $W_4$ are chosen to reflect the design specifications of robustness to parametric uncertainties, tracking bandwidth, and saturation limits on the control signal. More specifically, the weight functions $W_1(j\omega)$ and $W_2(j\omega)$ are chosen to be large in frequency range $\left[0,\omega_{BW}\right]$ to ensure small tracking errors $e_1 = V_{ref}-V_{dc}$ and $e_2 = i_{load} - i_L$ in this frequency range. The design of weight function $W_3(j\omega)$ entails ensuring that the control effort lies within saturation limits. The weight function $W_4$ is designed as a high-pass filter to ensure that the transfer function from $i_{load}$ to $V_{dc}$ is small at high frequencies to provide mitigation to measurement noise.

\begin{figure}
	\includegraphics[width= 0.95\columnwidth]{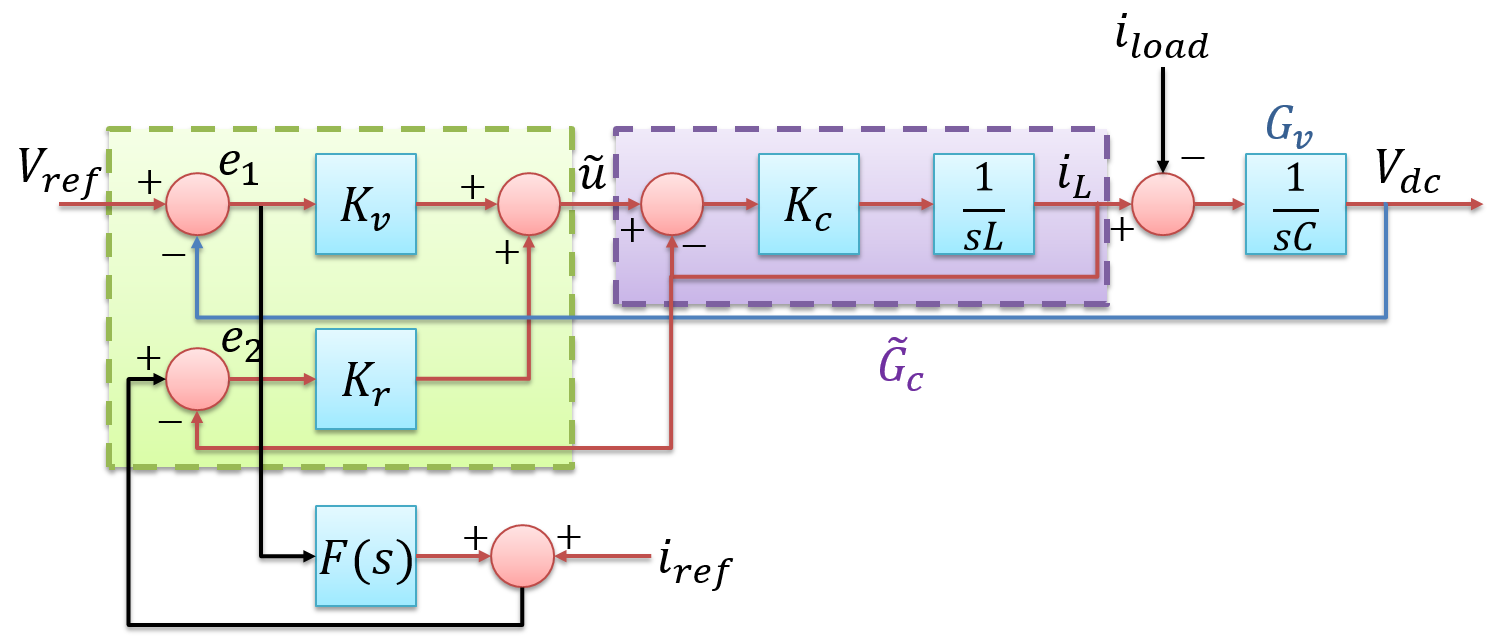}
	\caption{{\small Modified architecture for {\em decentralized} implementation. $F(s)$ represents the Laplace transform of the corresponding LTI filter $f(\cdot)$. The reference to the outer controller $K_r$ is $i_{ref} + f(\cdot)*(V_{ref}-V_{dc})$.}}
	\vspace{-1.5em}
	\label{fig:buck_mod}
\end{figure}

{\em Current droop compensation for voltage regulation without $i_{load}$ measurement}: So far in our analysis, we assume that $i_{load}$ is available for direct measurement for all converters. However in practice, while the converters can measure their common DC-link voltage $V_{dc}$, $i_{load}$ is only {\em estimated} through power calculations on the AC side and communicated to individual converters only at a rate slower than the sampling rate of the controllers. In the wake of this limitation, we propose a novel method for voltage regulation and power sharing that is inspired by conventional voltage droop method. An interesting aspect of the proposed implementation is that the same outer controllers $K_{v}$ and $K_{r}$ along with the shaped inner plant $\tilde{G}_{c,n}$ can be employed even for the scenario where $i_{load}$ measurement is unavailable. In this case, the outer controller $K_{r}$ regulates the inductor current $i_L$ to $i_{ref} + f(\cdot)*(V_{ref}-V_{dc})$, where $f(\cdot)$ represents an LTI filter and $*$ is the convolution operator. This can be understood as follows. Let us suppose that $i_{ref}>i_{load}=(V_{ref}/R)$, where $R$ is the output load resistance. This in turn implies that $V_{dc}=i_{ref}R$ is bigger than $V_{ref}$. Thus the compensation term $(V_{ref}-V_{dc})<0$ and the reference to the outer-loop controller $K_r$ becomes smaller than $i_{ref}$, as required. A similar inference can be drawn when $i_{ref}<i_{load}$. The corresponding block diagram is shown in Fig. \ref{fig:buck_mod}.

\begin{figure}
	\includegraphics[width= 0.95\columnwidth]{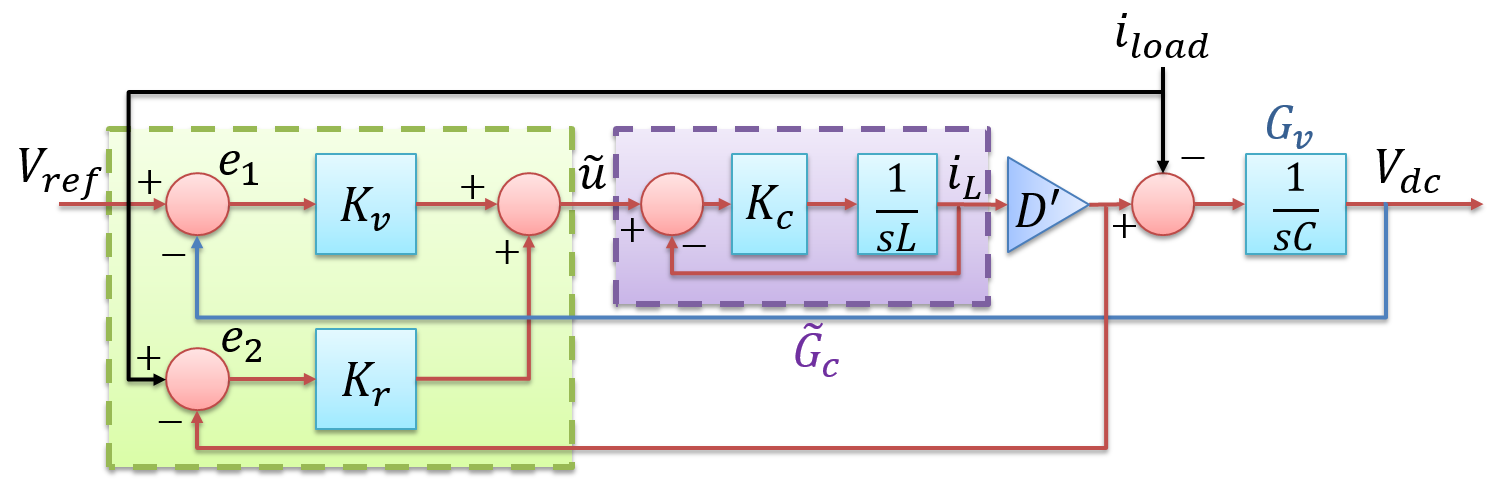}
	\caption{{\small Block diagram representation of the inner-outer control design for boost-converter. Note that this representation differs from that of the buck converter only in terms of the constant parameter $D'$ and the expression that is used to obtain the duty-cycle from the control signal $\tilde{u}$.}}
	\vspace{-1.5em}
	\label{fig:single_boost}
\end{figure}

{\em Extension to Boost Converters}: The extension of the proposed control design to Boost and Buck-Boost DC-DC converters is easily explained after noting that their averaged models are structurally identical to Buck converters, except that the dependence of duty cycles on the control signal $u$ or constant parameter $D'$ are different. The differences in how duty cycles depend on $u(t)$ do not matter from the control design viewpoint since duty cycles for pulse-width modulation are obtained only after obtaining the control designs (that use the averaged models). Fig. \ref{fig:single_boost} shows the equivalent schematic of the proposed control framework for a Boost converter system.

\section{EXTENSION TO A SYSTEM OF PARALLEL CONVERTERS}\label{sec:Control_Many}
In this section	we extend our control framework for a single converter to a system of DC-DC converters connected in parallel in the context of power sharing, keeping in mind the practicability and robustness to modeling and load uncertainties.
\begin{figure*}[!t]
	\begin{center}
	\begin{tabular}{cc}
	\includegraphics[width=0.9\columnwidth]{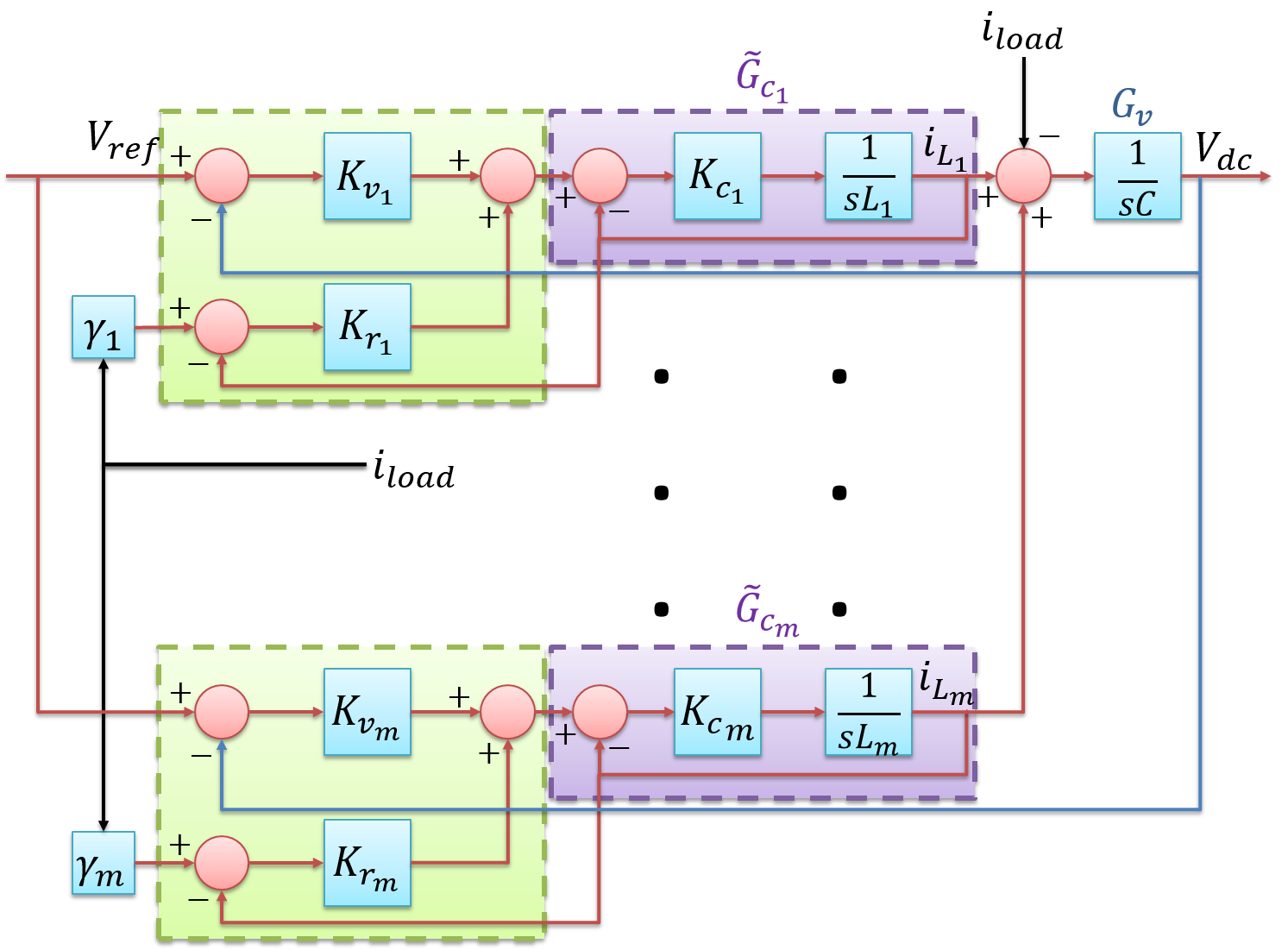}&\includegraphics[width=0.9\columnwidth]{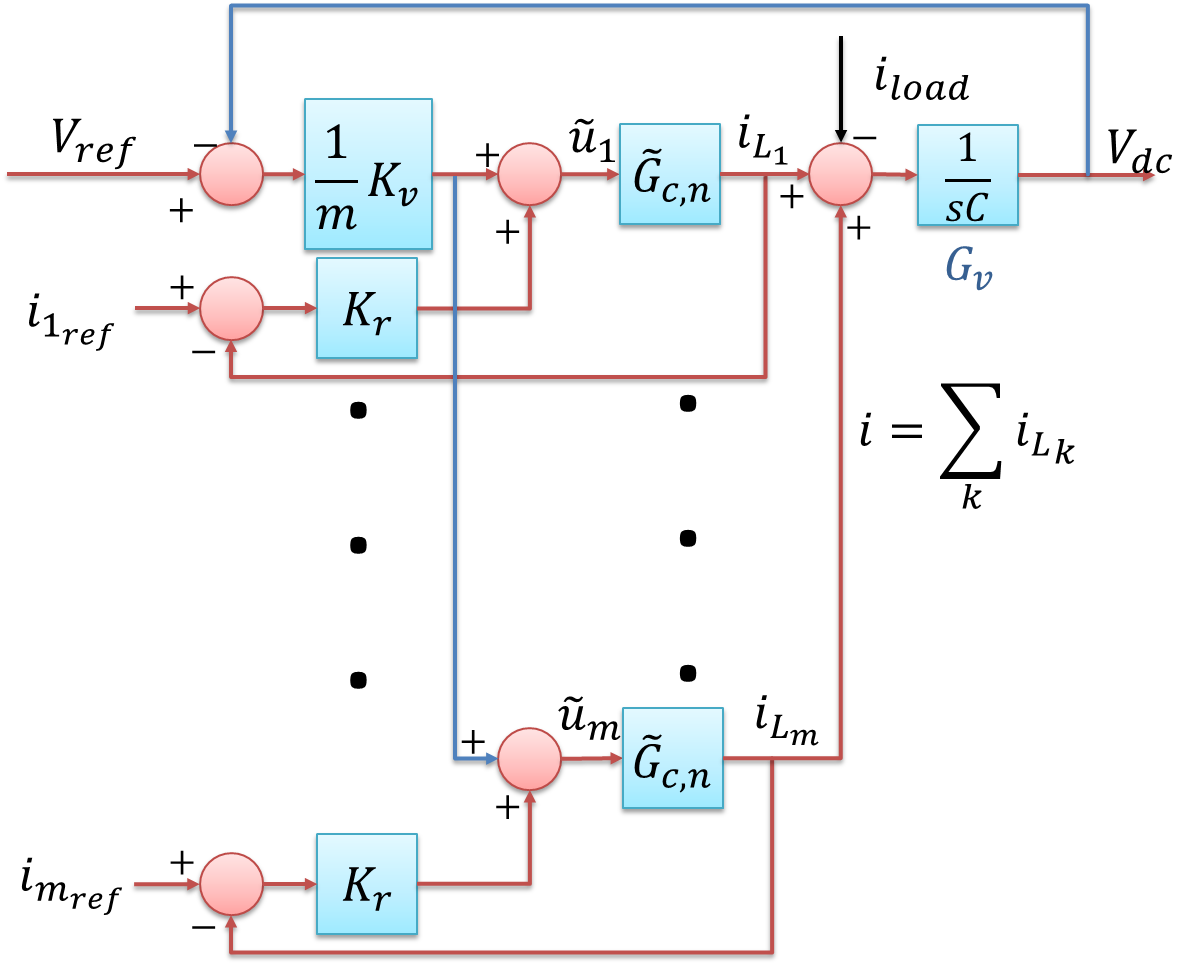}\cr
	(a)&(b)\end{tabular}
	\vspace{-0.5em}
	\caption{{\small (a) Control framework for a network of $m$ parallel converters. Here $\gamma_k$ represents the proportion of power demanded from the $k^{th}$ source. (b) A multiple-converters system with shaped inner plants $\tilde{G}_{c}$. In the proposed implementation, we adopt the same outer controller for different converters, i.e., $K_{v_1} = K_{v_2} = .. = K_{v_m} = \frac{1}{m}K_v$ and $K_{r_1} = K_{r_2} = .. = K_{r_m} = K_r$.}}
	\label{fig:model_many}
	\vspace{-2em}
	\end{center}
\end{figure*}

Fig. \ref{fig:model_many}a represents an inner-outer control framework for a system of $m$ parallel connected converters. Note that instead of feeding $i_{load}$ directly to the $k^{th}$ outer controller $K_{r_k}$, the measurement signal is prescaled by a time-varying multiplier $\gamma_k, 0\leq\gamma_k\leq 1$. The choice of $\gamma_k$ dictates the power sharing requirements on the $k^{th}$ converter. In fact, we later show that the proposed implementation distributes the output power in the ratios $\gamma_1:\gamma_2:..:\gamma_m$. After noting that the voltage-regulation and current reference tracking is common to all the outer controllers, in our architecture, we impose the same design for outer-controllers for all the converters, i.e., $K_{v_1} = K_{v_2} = .. = K_{v_m}$ and $K_{r_1} = K_{r_2} = .. = K_{r_m}$. This imposition enables significant reduction in the overall complexity of the distributed control design for a parallel network of converters and power sources, thus ensuring the practicability of the proposed design which allows integration of power sources of different types and values.

We design inner-controllers $K_{c_k}$ such that the inner-shaped plants from $\tilde{u}_k$ to $i_{L_k}$ are same and given by,
\begin{small}
	\begin{equation}\label{eq:Gcn}
		\tilde{G}_{c,n}(s) = \left(\frac{\tilde{\omega}}{s+\tilde{\omega}}\right)\left(\frac{s^2+2\zeta_{1,n}\omega_0s+\omega_0^2}{s^2+2\zeta_{2,n}\omega_0s+\omega_0^2}\right),
	\end{equation}
\end{small}
where the ratio $\zeta_{1,n}/\zeta_{2,n}$ determines the tradeoff of $120$Hz ripple between the total inductor current $i_L = \sum\limits_{k=1}^mi_{L_k}$ and the capacitor current $i_C$. Note that for given values of $\zeta_{1,n}, \zeta_{2,n}$ and inductance $L_k$, explicit design of $K_{c_k}$ exists and is given by (\ref{eq:Kc}). After noting that $K_{v_k} = \frac{1}{m}K_{v}$ and $K_{r_k} = K_{r}$, the system in Fig. \ref{fig:model_many}a can be simplified to Fig. \ref{fig:model_many}b.

Indeed, by our choice of inner and outer controllers, the transfer functions from external references $V_{ref}$ and $i_{k_{ref}}$ to the desired output $V_{dc}$ are identical for all converters. Hence the entire network of parallel converters can be analyzed in the context of an equivalent single converter system. This implies that $K_{v_k}$ and $K_{r_k}$ can be computed by solving $\mathcal{H}_\infty$-optimization problem (as discussed in the previous section) similar to the {\em single} converter case. We make these design specifications more precise and bring out the equivalence of the control design for the single and multiple converter systems in the following theorem.

We say that the system representation in Fig. \ref{fig:single_block} is {\em equivalent} to that in Fig. \ref{fig:model_many}b, when the transfer functions from the reference voltage $V_{ref}$ and load current $i_{load}$ to the DC-link voltage $V_{dc}$ (and therefore the total current sourced $i=i_L$) in Fig. \ref{fig:single_block} are identical to the corresponding transfer functions in Fig. \ref{fig:model_many}b. In the following theorem, we denote the outer-voltage sensitivity transfer function in Fig. \ref{fig:single_block} from $V_{ref}$ to $V_{dc}$ by $S_1=\left(\frac{1}{1\!+\!\tilde{G}_{c,n}K_r\!+\!G_v\tilde{G}_{c,n}K_v}\right)$ and the corresponding complementary sensitivity transfer function by $T_1=G_v\tilde{G}_{c,n}K_vS_1$. Similarly, we denote the outer-current sensitivity transfer function by $S_2=\left(\frac{1}{1\!+\!\tilde{G}_{c,n}K_r}\right)$ and the corresponding complementary transfer function by $T_2=1-S_2$. Moreover, we define $H=\tilde{G}_{c,n}K_vS_1$.

\begin{theorem}\label{thm1}
Consider the single-converter system in Figure \ref{fig:single_block} with inner-shaped plant $\tilde{G}_{c,n}(s)$ as given in (\ref{eq:Gcn}), outer controllers $K_{v}$, $K_{r}$ and external references $V_{ref}$, $i_{load}$; and the multi-converter system described in Figures \ref{fig:model_many}a and \ref{fig:model_many}b with inner-shaped plants $\tilde{G}_{c_k}=\tilde{G}_{c,n}(s)$ and outer controllers $K_{v_k}=\frac{1}{m}K_v$; $K_{r_k}=K_{r}$, and external references $V_{ref}$, $i_{k_{ref}}$ for $1\leq k\leq m$.\\
\noindent 1. [System Equivalence]: If $\sum_{k=1}^{m}i_{k_{ref}}=i_{load}$, then the system representations in Fig. \ref{fig:single_block} and Fig. \ref{fig:model_many}b are {\em equivalent}.\\
\noindent 2. [Power Sharing]: If controllers $K_v$ and $K_r$ are designed such that $|H(j\omega)|< \epsilon$, $|S_2(j\omega)|< \epsilon$ and the total current mismatch $|\sum_{k}i_{k_{ref}}(j\omega)-i_{load}(j\omega)|<\Delta$ at frequncy $\omega$ for some $\epsilon>0$ and $\Delta>0$, then $|i_{L_k}(j\omega)-i_{k_{ref}}(j\omega)|< \frac{\epsilon(|V_{ref}(j\omega)|+m|i_{k_{ref}}(j\omega)|+(1+\epsilon)|i_{load}(j\omega)|)}{m}+\frac{(1+\epsilon)^2\Delta}{m}$.
\end{theorem}
{\bf Remark 1:} Since the system representations in Figs. \ref{fig:single_block} and \ref{fig:model_many}b are equivalent, the analysis and design of the entire multi-converter system can be done using an equivalent single converter system, where the multi-converter system inherits the performance and robustness achieved by a design for the single-converter system.\\
{\bf Remark 2:} While the sensitivity transfer functions $S_1$ and $S_2$ can be made sufficiently small by designing appropriate high DC-gain controllers, the transfer function $H$ is also small at sufficiently low-frequencies. In fact, it can be easily shown that $|H(j0)|=0$, since $|\tilde{G}_{c,n}(j0)|=1$ and $G_v(s)=1/(sC)$, which has infinite DC-gain.\\
{\bf Remark 3:} Note that from power-sharing result in the above theorem, if $i_{k_{ref}}=\gamma_ki_{load}$, where $\sum_{k}\gamma_k=1, \gamma_k\geq 0$, then the output current at the DC-link gets divided approximately in the ratio $\gamma_1:\gamma_2:\dots :\gamma_m$; more precisely the {\em low}-frequency components (and thus the steady-state) $i_{L_1}:i_{L_2}:\dots :i_{L_m} \approx \gamma_1:\gamma_2:\dots :\gamma_m$.\\
Proof: See Appendix.

\section{CASE STUDIES: SIMULATIONS AND DISCUSSIONS}\label{sec:simulations}
In this section, we report some simulation studies that cover different aspects of the proposed distributed control design. All simulations are performed in MATLAB/Simulink using SimPower/SimElectronics library. Note that the experiments are underway and therefore not reported in this paper. In order to include nonlinearities associated with real-world experiments, and effects of switching frequencies on voltage regulation and power sharing, we use {\em non-ideal} components (such as diodes with non-zero forward-bias voltage, IGBT switches, stray capacitances, parametric uncertainties) and switched level implementation.

\subsection{Robustness to Modeling Parameters}\label{subsec:simulations_robustness}
Traditional control techniques such as proportional-integral (PI) based control designs exhibit satisfactory performance when the  actual converter system parameters ($L, C$) lie {\em close} to the nominal system parameters for which the controllers are to be designed. A slight deviation from the nominal values may result in rapid degradation in the tracking performance and power sharing. The issue is particularly critical in power electronic systems where individual component values have large tolerance about the nominal values. The lack of robustness is addressed through $\mathcal{H}_\infty$ robust control framework, where an optimizing controller with guaranteed margins of robustness to modeling uncertainties is sought.

\begin{figure*}[!t]
	\begin{center}
	\begin{tabular}{cc}
	\includegraphics[width=0.85\columnwidth]{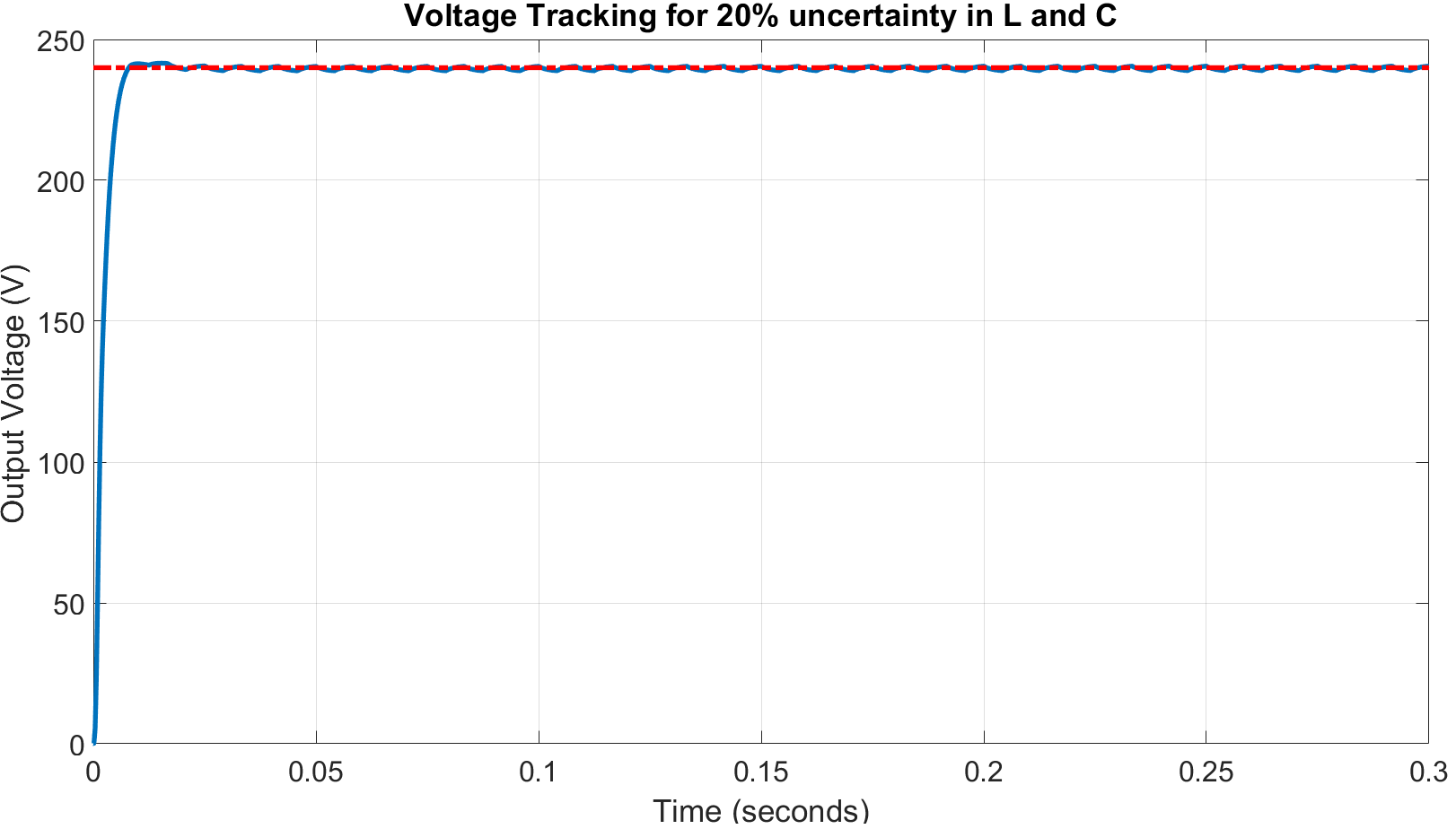}&\includegraphics[width=0.9\columnwidth]{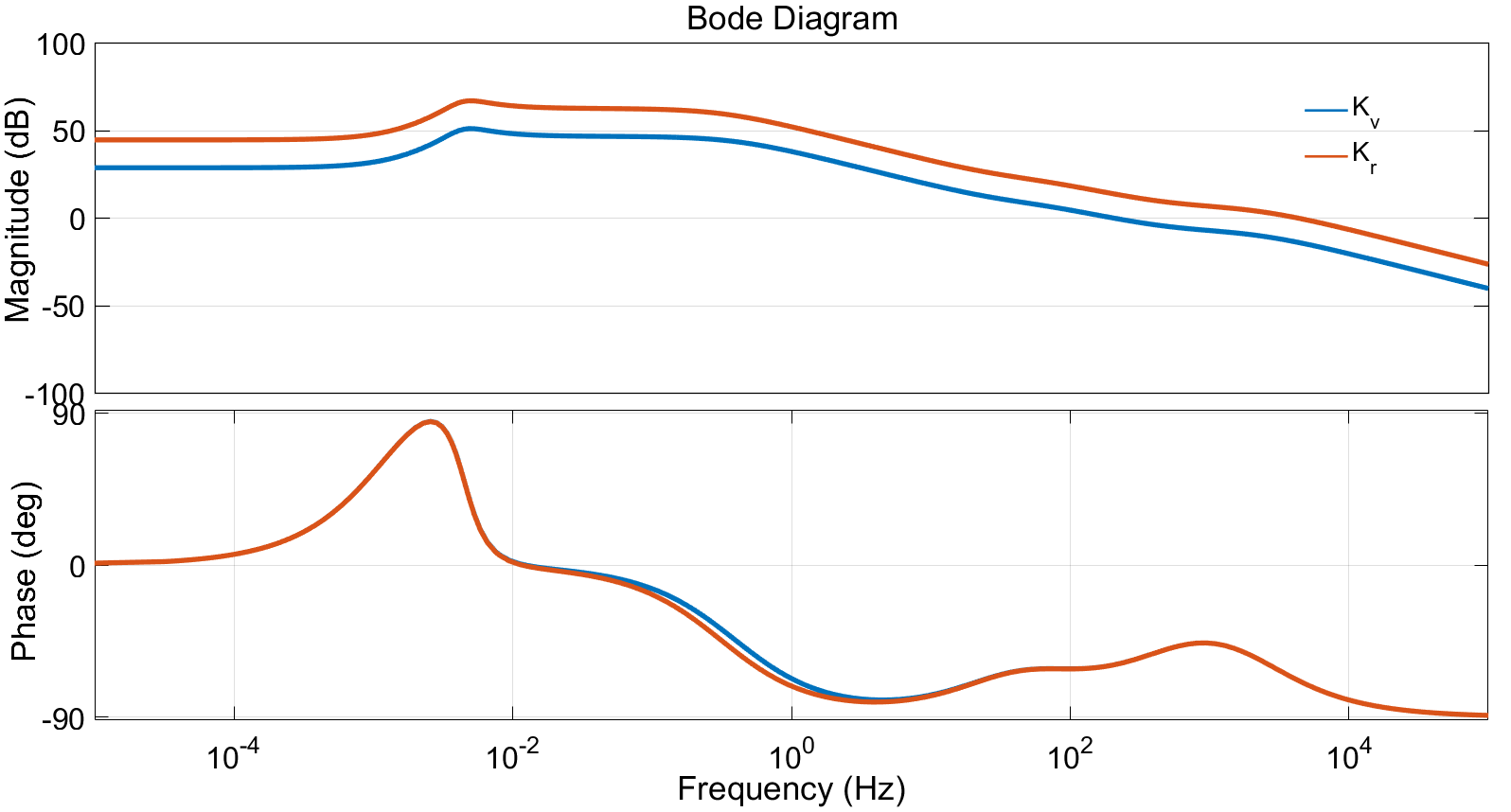}\cr
	(a)&(b)\end{tabular}
	\vspace{-0.5em}
	\caption{{\small (a) Robustness to model parameters. The proposed controller tracks the desired voltage in presence of parametric uncertainties and ripple in load current. (b) Bode plots of the outer-loop controllers. It can be observed that the two controller transfer functions are constant multiples of one-another.}}
	\label{fig:sim_uncertain}
	\vspace{-2em}
	\end{center}
\end{figure*}

Fig. \ref{fig:sim_uncertain}a shows the tracking performance of the proposed robust inner-outer controllers for $20\%$ uncertainty in indcutance ($L$) and capacitance ($C$) values. The actual converter parameters are chosen as:\\
$C = 500 \mu$F, $L = 1.2$mH, Switching-frequency $f_s = 50$kHz, Input voltage $V_g = 480V$, Desired output voltage $V_{ref} = 240V$, Load-current $i_{load} = 20 + 0.4\sin(2\pi 120t)$\\
The design parameters for the inner-controller $K_c$ are:\\
Damping factors $\zeta_1 = 1.2$, $\zeta_2 = 2.1$, and $\tilde{\omega} = 2\pi 200$rad/s.\\
The outer controllers $K_{v}$ and $K_{r}$ are obtained by solving the stacked $\mathcal{H}_\infty$ optimization problem (see Eq. (\ref{eq:hinf}))\cite{skogestad2007multivariable} with the weighting functions:
\begin{small}
	\begin{equation*}
		W_1 = \frac{0.5(s+502.7)}{(s+2.513)}, W_2 = \frac{0.5(s+628.3)}{(s+3.142)}, W_3 = 0.1
	\end{equation*}
\end{small}
The resulting outer-controllers are reduced to sixth-order using balanced reduction \cite{dullerud2013course} and are given by,
\begin{small}
	\begin{eqnarray}
		K_{v} &=& -0.0076\frac{(s-8.69e5)}{(s+2.73e4)}\frac{(s+2.01e4)}{(s+1.07e4)}\frac{(s+2577)}{(s+433.9)}\nonumber \\
		&&\frac{(s+194.2)}{(s+2.498)}\frac{(s^2+0.02s+0.0001)}{(s^2+0.01978s+0.0008)}\nonumber\\
		K_{r} &=& 0.065\frac{(s+4.07e5)}{(s+1.15e4)}\frac{(s+2474)}{(s+422.4)}\frac{(s+191.7)}{(s+3.11)} \nonumber \\
		&&\frac{(s+3.20)}{(s+2.03)}\frac{(s+0.01)(s+0.0099)}{(s^2+0.01978s+0.0008)}\nonumber 
	\end{eqnarray}
\end{small}

{\em Optimal structure of outer-loop controllers}: For the stacked $\mathcal{H}_\infty$ problem, the outer-loop controllers $K_{v}$ and $K_{r}$ are observed to be constant multiples of each other for a number of choices of weighting functions $W_1$ and $W_2$. While we are yet to explore the reasons for the underlying optimal structure, the optimal structure significantly reduces the complexity of the distributed control design by allowing to get rid of the controller $K_{r}$ from the outer-loop and modifying the input of the controller $K_{v}$ to $V_{ref}-V_{dc} + \alpha(i_{load}-i_L)$. Here $\alpha > 0$ is an appropriate constant which captures the relationship between the outer-loop controllers $K_{v}$ and $K_{r}$. We believe that the constant $\alpha$ is related to the system parameters $L$ and $C$, however, this is something we would definitely like to explore in our future work.

\subsection{Current Sharing among Converters}\label{subsec:simulations_sharing}
\begin{figure}
	\includegraphics[width= 0.95\columnwidth]{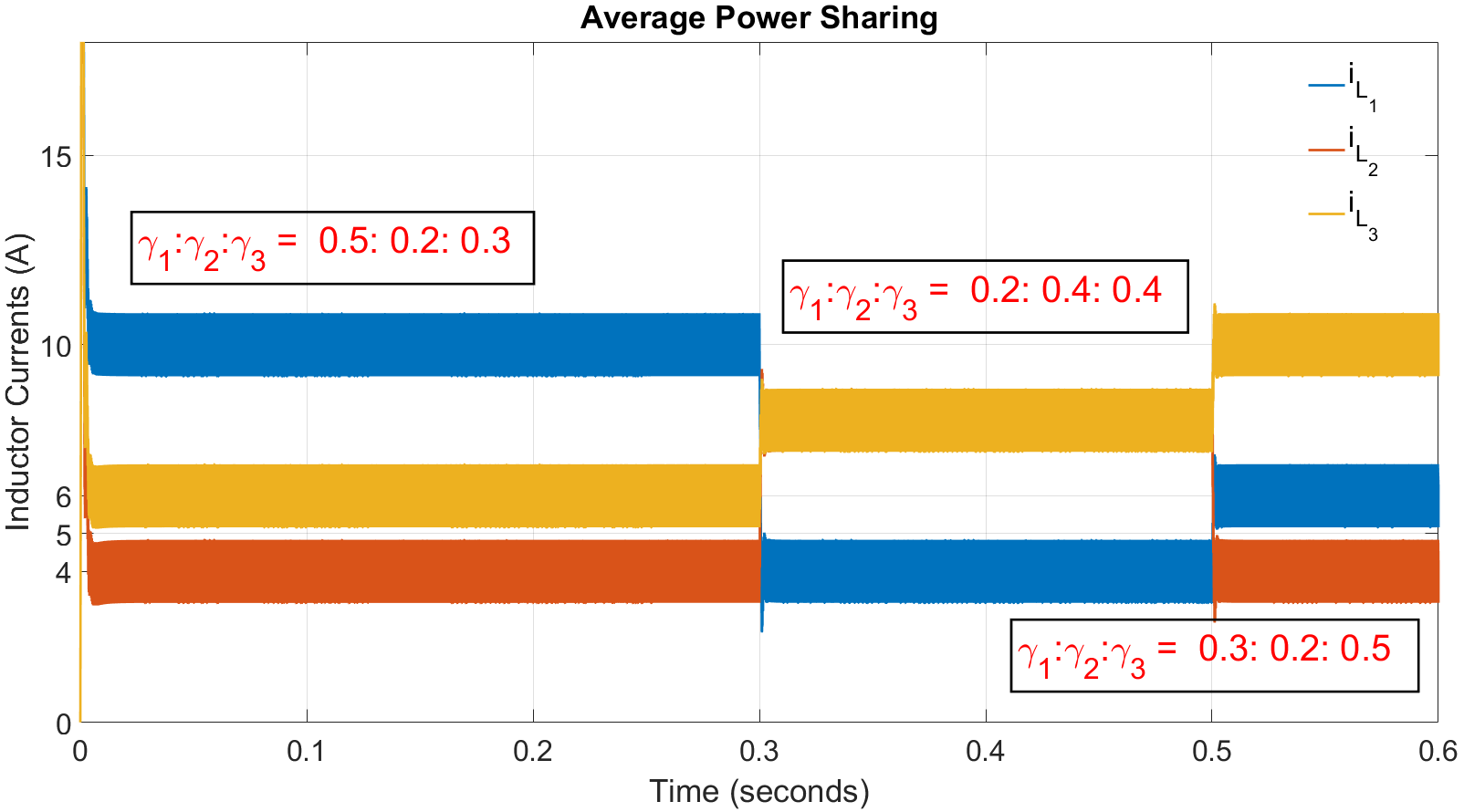}
	\caption{{\small Power sharing among three buck-converters for $V_{ref} = 240V$ and $R = 12\Omega$ in the prescribed ratios. The controllers allow for rapidly varying power-sharing requirements.}}
	\vspace{-1.5em}
	\label{fig:curr_sharing}
\end{figure}

Fig. \ref{fig:curr_sharing} shows the system performance for time-varying output current sharing among three Buck converters. The following parameters are assumed for the three converters - $L_1 = 1.2$mH, $V_{g_1} = 480V$; $L_2 = 1.6$mH, $V_{g_2} = 460V$; $L_3 = 1.9$mH, $V_{g_3} = 480V$.
Clearly, the total load current is $20$A. The converters divide the load current in the ratios $10:4:6$ for $t\in[0,0.3s]$; $4:8:8$ for $t\in[0.3s,0.5s]$ and $6:4:10$ for $t\in[0.5s,0.6s]$, as required.

\subsection{Current Sharing and Voltage Regulation without $i_{load}$ measurement}\label{subsec:simulations_no_iload}
\begin{figure}
	\includegraphics[width= 0.95\columnwidth]{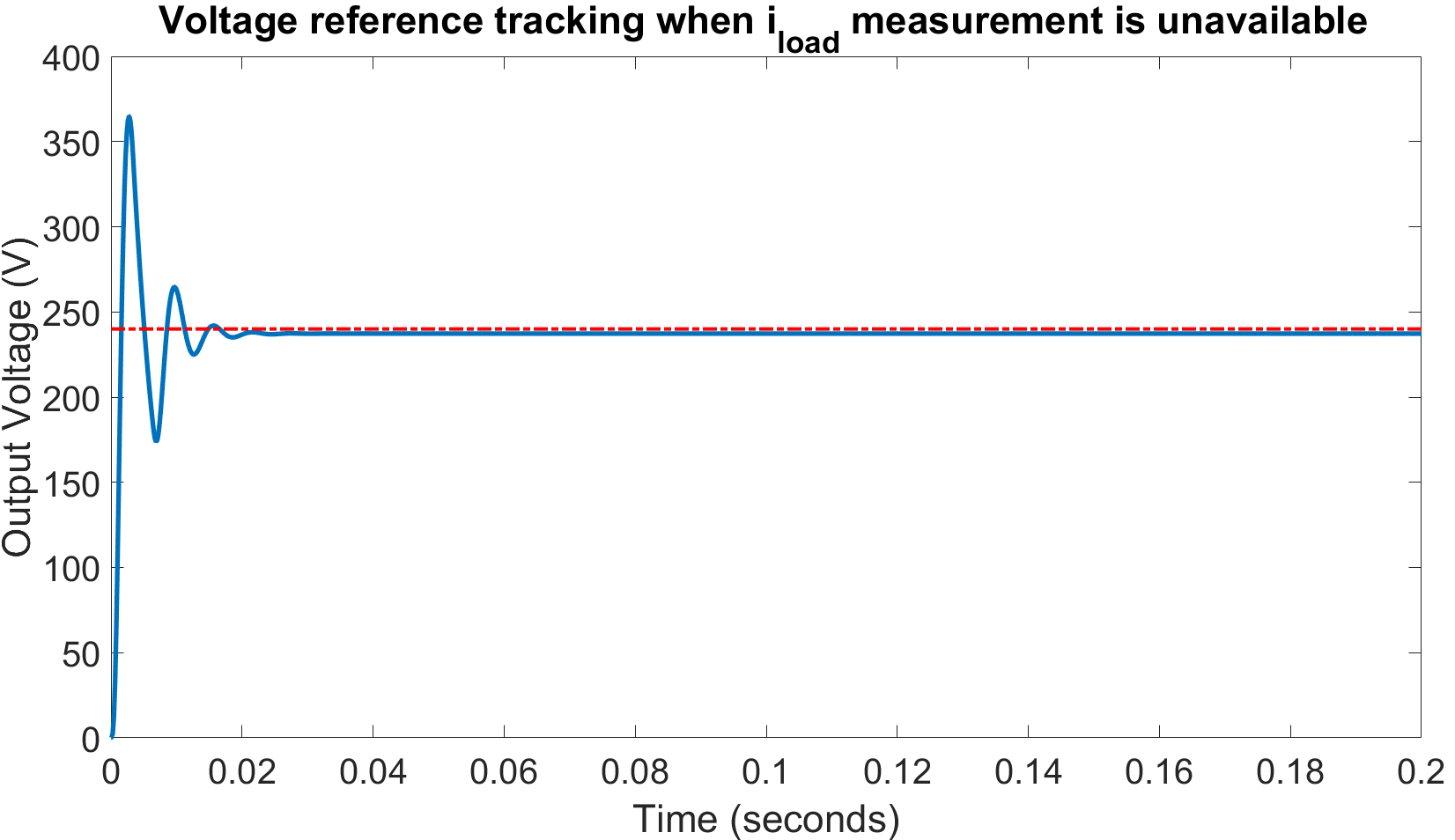}
	\caption{{\small A droop-like approach for voltage regulation in the absence of $i_{load}$ measurement. The converter is required to regulate the output voltage to $240V$. The load $R$ is chosen to be $12\Omega$. This corresponds to an $i_{load} = 20A$. However, in the absence of $i_{load}$ measurement, we assume $i_{ref} = 16A$. The filter transfer function is chosen as $F(s) = \frac{376.99}{s+314.16}$.}}
	\vspace{-1.5em}
	\label{fig:no_iload}
\end{figure}

We now consider the scenario when $i_{load}$ measurement is unavailable. In this case, the input to the outer-loop controller $K_{r}$ is some nominal reference current  $i_{ref}$ plus a compensation term which is a manifestation of error in voltage reference tracking. This can be understood as follows. Fig. \ref{fig:no_iload} shows the system performance for voltage reference tracking for the proposed {\em decentralized} implementation. Note that despite the unavailability of load current measurement, the controller regulates the DC-link voltage to the desired reference voltage, albeit with relatively larger overshoot.

\section{CONCLUSIONS AND FUTURE WORKS}\label{sec:conclusion}
In this work, we propose a distributed control architecture for voltage tracking and power sharing for a network of DC-DC converters connected in parallel. The proposed design is capable of achieving multiple objectives such as robustness to modeling uncertainties, reference DC voltage generation and output power sharing among multiple DC sources. The controllers are designed using a robust optimal control framework. We also propose a novel approach for {\em decentralized} implementation, where the load current is not available for measurement. We are currently setting up the experiments to demonstrate the effectiveness of the proposed implementation. Moreover, the optimal structure of the outer-loop controllers $K_{v}$ and $K_{r}$ needs to be analyzed in full details to gain further insights into the problem of voltage control of DC-DC converters.

\section*{APPENDIX}
\subsection*{Proof of Theorem 1: System Equivalence}
\begin{proof}
	The underlying equivalence is straightforward to derive. From Fig. \ref{fig:single_block}, with $\tilde{G}_c(s) = \tilde{G}_{c,n}(s)$, the output DC-link voltage in terms of the exogenous signals $V_{ref}$ and $i_{load}$ is given by
	\begin{small}
		\begin{equation}\label{eq:app1}
			V_{dc} = \left(\frac{1}{1\!+\!\tilde{G}_{c,n}K_r\!+\!G_v\tilde{G}_{c,n}K_v}\right)\left(G_v\tilde{G}_{c,n}K_vV_{ref}-G_vi_{load}\right)
		\end{equation}
	\end{small}
	However, from Figs. \ref{fig:model_many}a and \ref{fig:model_many}b, we obtain
	\begin{small}
		\begin{eqnarray}\label{eq:app2}
			V_{dc} = G_v\left(-i_{load}+\tilde{G}_{c,n}\sum_{k=1}^{m}\tilde{u}_k\right)\nonumber \\
			\tilde{u}_k = \left(\frac{1}{1+\tilde{G}_{c,n}K_r}\right)\left(\frac{1}{m}K_v(V_{ref}-V_{dc})+K_ri_{k_{ref}}\right)
		\end{eqnarray}
	\end{small}
	From Eq. (\ref{eq:app2}) and using the fact that $\sum_{k=1}^{m}i_{k_{ref}}=i_{load}$ we recover Eq. (\ref{eq:app1}), which establishes the required equivalence.
\end{proof}

\subsection*{Proof of Theorem 1: Power Sharing}
\begin{proof}
	From Fig. \ref{fig:sim_uncertain}b, it can be shown that the difference between the inductor current $i_{L_k}$ and the reference current $i_{k_{ref}}$ for the $k^{th}$ converter (in terms of signals $V_{ref}$, $i_{load}$ and $i_{k_{ref}}$) is given by
	\begin{small}
		\begin{eqnarray}\label{eq:app3}
			i_{L_k}-i_{k_{ref}}=&&\frac{1}{m}HV_{ref}-\frac{1}{m}T_1T_2\left(\sum_ki_{k_ref}-i_{load}\right) \nonumber \\
			&&+ \frac{1}{m}T_1S_2i_{load} - S_2i_{k_ref}
		\end{eqnarray}
	\end{small}
	Moreover, we have $|T_1(j\omega)|< 1+\epsilon$ and $|T_2(j\omega)|< 1+\epsilon$. Thus from (\ref{eq:app3}) and conditions of the theorem, we get
	\begin{small}
		\begin{eqnarray}\label{eq:app4}
			|i_{L_k}(j\omega)-i_{k_{ref}}(j\omega)|< \frac{\epsilon}{m}|V_{ref}(j\omega)| + \frac{\epsilon(1+\epsilon)}{m}|i_{load}(j\omega)| \nonumber \\
			+ \frac{(1+\epsilon)^2\Delta}{m} + \epsilon |i_{k_{ref}}(j\omega)|
		\end{eqnarray}
	\end{small}
	Additionally, if $i_{k_{ref}}=\gamma_ki_{load}$, then the total current mismatch term in (\ref{eq:app3}) is zero. Moreover, through an appropriate design of high DC-gain controllers, we have that there exists a sufficiently small $\epsilon >0$ such that for $\omega << \omega_{BW}$, where $\omega_{BW}$ is bandwidth of the closed-loop system, $|H(j\omega)|<\epsilon$ and  $|S_2(j\omega)|<\epsilon$. Then from (\ref{eq:app4}) we obtain,
	\begin{small}
		\begin{eqnarray}
			|i_{L_k}(j\omega)\!-\!\gamma_ki_{load}(j\omega)|< \frac{\epsilon}{m}|V_{ref}(j\omega)|\!+\!\frac{\epsilon}{m}|i_{load}(j\omega)|\nonumber \\
			+\epsilon |i_{k_{ref}}(j\omega)|
		\end{eqnarray}
	\end{small}
	Thus $|i_{L_k}(j\omega)\!-\!\gamma_ki_{load}(j\omega)|$ is small for $\omega<<\omega_{BW}$ and therefore, the load current gets approximately divided in the ratio $\gamma_1:\gamma_2:\dots :\gamma_m$.
\end{proof}

\bibliographystyle{IEEEtran} 
\bibliography{myRef}

\end{document}